\documentclass[a4paper]{llncs}

\usepackage[latin1]{inputenc}
\usepackage[T1]{fontenc}
\usepackage{ae,aecompl}
\usepackage{amsmath,amsfonts,amssymb}
\usepackage[tableposition=bottom,font=small,labelfont=bf,format=hang]{caption}
\usepackage{booktabs}
\usepackage{subfig}
\usepackage[dvips]{graphicx}
\usepackage[dvips,colorlinks=false,bookmarks=false]{hyperref}

\hypersetup{
pdfauthor = {Lars Eirik Danielsen and Matthew G. Parker},
pdftitle = {Edge Local Complementation and Equivalence of Binary Linear Codes},
pdfkeywords = {binary linear codes, graphs, edge local complementation, classification}
}

\DeclareMathOperator{\GF}{GF}
\DeclareMathOperator{\wt}{wt}

\begin{document}

\pagestyle{headings}

\mainmatter

\title{Edge Local Complementation and Equivalence of Binary Linear Codes}
\author{Lars Eirik Danielsen \and Matthew G. Parker}
\institute{Department of Informatics, University of Bergen, PB 7803, N-5020 Bergen, Norway\\
\texttt{\{\href{mailto:larsed@ii.uib.no}{larsed},\href{mailto:matthew@ii.uib.no}{matthew}\}@ii.uib.no}\\
\texttt{http://www.ii.uib.no/\~{}\{\href{http://www.ii.uib.no/~larsed}{larsed},\href{http://www.ii.uib.no/~matthew}{matthew}\}}}

\maketitle

\begin{abstract}
Orbits of graphs under the operation \emph{edge local complementation} (ELC)
are defined. We show that the ELC orbit of a \emph{bipartite} graph corresponds
to the equivalence class of a \emph{binary linear code}. The \emph{information sets} 
and the \emph{minimum distance} of a code can
be derived from the corresponding ELC orbit. By extending earlier results 
on \emph{local complementation} (LC) orbits, we classify the ELC orbits of all graphs 
on up to 12 vertices. We also give a new method for classifying binary linear codes, 
with running time comparable to the best known algorithm.\\[\baselineskip]
\end{abstract}

\section{Introduction}

In this section we first give some definitions from graph theory, in particular 
we describe the two graph operations \emph{local complementation} (LC) and 
\emph{edge local complementation} (ELC), the latter also known as the \emph{pivot}
operation.
We then give some definitions related to \emph{binary linear codes}. Of particular
interest is the concept of \emph{code equivalence}.
Östergård~\cite{ostergard} represented codes as graphs, and devised an algorithm
for classifying codes up to equivalence.
In Section~\ref{sec:equiv}, we show a different way of representing a binary linear code as
a \emph{bipartite} graph. We prove that ELC on this graph provides a simple
way of jumping between equivalent codes, and that the orbit of a bipartite graph under
ELC corresponds to the complete equivalence class of the corresponding code.
We also show how ELC on a bipartite graph generates all \emph{information sets} 
of the corresponding code. Finally, we show that the \emph{minimum distance} of
a code is related to the minimum vertex degree over the corresponding ELC orbit.
In Section~\ref{sec:class} we describe our algorithm for classifying ELC orbits,
which we have used to generate all ELC orbits of graphs on up to 12 vertices.
Although ELC orbits of non-bipartite graphs do not have any obvious applications to
classical coding theory, they are of interest in other contexts, such as
\emph{interlace polynomials}~\cite{interlace2,interlace} 
and \emph{quantum graph states}~\cite{nestELC} which are related to \emph{quantum error correcting codes}.
From the ELC orbits of bipartite graphs a classification of binary linear codes 
can be derived. Binary linear codes have previously been classified up to 
length 14~\cite{ostergard,ostergard2}. We have generated the bipartite ELC orbits of 
graphs on up to 14 vertices, and this classification can be extended to at least 15 
vertices [Sang-il Oum, personal communication], showing that our method
is comparable to the best known algorithm.
However, the main result of this paper is not a classification of codes, but
a new way of representing equivalence classes of codes, and a classification of
all ELC orbits of length up to 12.

\subsection{Graph Theory}

A \emph{graph} is a pair $G=(V,E)$ where $V$ is a set of \emph{vertices},
and $E \subseteq V \times V$ is a set of \emph{edges}. A graph with $n$ vertices
can be represented by an $n \times n$ \emph{adjacency matrix} $\Gamma$, where
$\gamma_{ij} = 1$ if $\{i,j\} \in E$, and $\gamma_{ij} = 0$ otherwise.
We will only consider \emph{simple} \emph{undirected} graphs whose
adjacency matrices are symmetric with all diagonal elements being 0, i.e., 
all edges are bidirectional and no vertex can be adjacent to itself.
The \emph{neighbourhood} of $v \in V$, denoted $N_v \subset V$, is the set of 
vertices connected to $v$ by an edge. 
The number of vertices adjacent to $v$ is called the \emph{degree} of $v$.
The \emph{induced subgraph} of $G$ on $W \subseteq V$ 
contains vertices $W$ and all edges from $E$ whose endpoints are both in $W$.
The \emph{complement} of $G$ is found by replacing $E$ with $V \times V - E$,
i.e., the edges in $E$ are changed to non-edges, and the non-edges to edges.
Two graphs $G=(V,E)$ and $G'=(V,E')$ are \emph{isomorphic} if and only if
there exists a permutation $\pi$ on $V$ such that $\{u,v\} \in E$ if and 
only if $\{\pi(u), \pi(v)\} \in E'$.
A \emph{path} is a sequence of vertices, $(v_1,v_2,\ldots,v_i)$, such that
$\{v_1,v_2\}, \{v_2,v_3\},$ $\ldots, \{v_{i-1},v_{i}\} \in E$.
A graph is \emph{connected} if there is a path from any vertex to any other vertex in the graph.
A graph is \emph{bipartite} if its set of vertices can be decomposed into two disjoint sets 
such that no two vertices within the same set are adjacent. We call a graph
\emph{$(a,b)$-bipartite} if its vertices can be decomposed into sets of size $a$ and $b$.

\begin{figure}
 \centering
 \subfloat[The Graph $G$]
 {\hspace{5pt}\includegraphics[width=.25\linewidth]{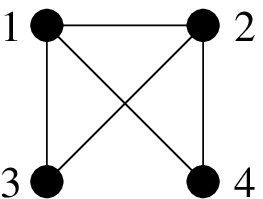}\hspace{5pt}\label{fig:lcexample1}}
 \quad
 \subfloat[The Graph $G*1$]
 {\hspace{5pt}\includegraphics[width=.25\linewidth]{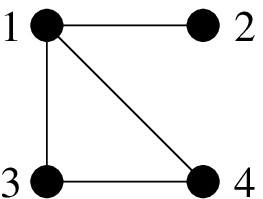}\hspace{5pt}\label{fig:lcexample2}}
 \caption{Example of Local Complementation}\label{fig:lcexample}
\end{figure}

\begin{definition}[\hspace{1pt}\hspace{-1pt}\cite{flaas,bouchet,hubert}]
Given a graph $G=(V,E)$ and a vertex $v \in V$,
let $N_v \subset V$ be the neighbourhood of $v$.
\emph{Local complementation} (LC) on $v$ transforms $G$ into $G * v$ by
replacing the induced subgraph of $G$ on $N_v$ by its complement. (Fig.~\ref{fig:lcexample})
\end{definition}

\begin{definition}[\hspace{1pt}\hspace{-1pt}\cite{bouchet}]\label{prop:triplelc}
Given a graph $G=(V,E)$ and an edge $\{u,v\} \in E$, \emph{edge local complementation} (ELC)
on $\{u,v\}$ transforms $G$ into $G^{(uv)} = G*u*v*u = G*v*u*v$.
\end{definition}

\begin{figure}
 \centering
 \includegraphics[width=.40\linewidth]{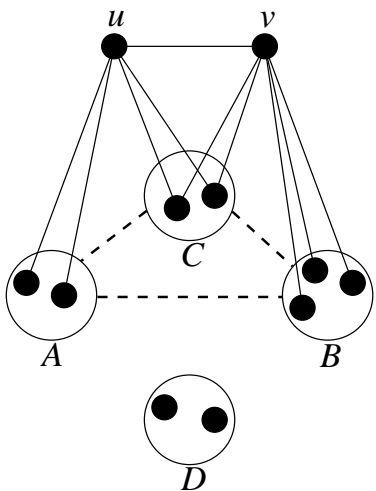}
 \caption{Visualization of the ELC Operation}\label{fig:elc}
\end{figure}

\begin{definition}[\hspace{1pt}\hspace{-1pt}\cite{bouchet}]\label{def:elc}
ELC on $\{u,v\}$ can equivalently be defined as follows. Decompose $V\setminus \{u,v\}$ 
into the following four disjoint sets, as visualized in Fig.~\ref{fig:elc}.
\renewcommand{\labelenumi}{$\Alph{enumi}$}
\begin{enumerate}
\item Vertices adjacent to $u$, but not to $v$.
\item Vertices adjacent to $v$, but not to $u$.
\item Vertices adjacent to both $u$ and $v$.
\item Vertices adjacent to neither $u$ nor $v$.
\end{enumerate}
To obtain $G^{(uv)}$, perform the following procedure.
For any pair of vertices $\{x,y\}$, where $x$ belongs to class $A$, $B$, or $C$,
and $y$ belongs to a different class $A$, $B$, or $C$, ``toggle'' the pair $\{x, y\}$,
i.e., if $\{x,y\} \in E$, delete the edge, and if $\{x,y\} \not\in E$, add the edge
$\{x,y\}$ to $E$. Finally, swap the labels of vertices $u$ and $v$.
\end{definition}

\begin{definition}
The \emph{LC orbit} of a graph $G$ is the set of all graphs that
can be obtained by performing any sequence of LC operations on~$G$.
Similarly, the \emph{ELC orbit} of $G$ comprises all graphs that
can be obtained by performing any sequence of ELC operations on~$G$.
(Usually we consider LC and ELC orbits of unlabeled graphs. In the cases
where we consider orbits of labeled graphs, this will be noted.)
\end{definition}

The LC operation was first defined by de Fraysseix~\cite{hubert}, and later studied
by Fon-der-Flaas~\cite{flaas} and Bouchet~\cite{bouchet}. Bouchet defined ELC 
as ``complementation along an edge''~\cite{bouchet}, but this operation is also 
known as \emph{pivoting} on a graph~\cite{interlace2,riera}.
LC orbits of graphs have been used to study \emph{quantum graph states}~\cite{glynnbook,hein,nest2},
which are equivalent to \emph{self-dual additive codes over $\GF(4)$}~\cite{calderbank}.
We have previously used LC orbits to classify such codes~\cite{setapaper,selfdual}.
ELC~orbits have also been studied in the context of quantum graph states~\cite{nestELC,riera}.
\emph{Interlace polynomials} of graphs have been defined with respect to both ELC~\cite{interlace2}
and LC~\cite{interlace}. These polynomials encode properties of the graph orbits, and 
were originally used to study a problem related to DNA sequencing~\cite{interlacedna}.

\begin{proposition}
If $G=(V,E)$ is a connected graph, then, for any vertex $v \in V$, $G*v$ must also be connected.
Likewise, for any edge $\{u,v\} \in E$, $G^{(uv)}$ must be connected.
\end{proposition}
\begin{proof}
If the edge $\{x,y\}$ is deleted as part of an LC operation on $v$, both $x$ and $y$ must be, and 
will remain, connected to $v$. 
Similarly, if by performing ELC on the edge $\{u,v\}$, the edge $\{x,y\}$ is deleted,
both $x$ and $y$ will remain connected to either $u$, $v$, or both, and $u$ and $v$ will remain
connected.\qed
\end{proof}

\begin{proposition}[\hspace{1pt}\hspace{-1pt}\cite{riera}]\label{prop:bipartite}
If $G$ is an $(a,b)$-bipartite graph, then, for any edge $\{u,v\} \in E$, $G^{(uv)}$ must also be 
$(a,b)$-bipartite.
\end{proposition}
\begin{proof}
A bipartite graph with an 
edge $\{u,v\}$ can not contain any vertex that is connected to both $u$ and $v$.
Using the terminology of Definition~\ref{def:elc}, the set $C$ will always be empty when
we perform ELC on a bipartite graph.
Moreover, all vertices in the set $A$ must belong to the same partition as $u$,
and all vertices in $B$ must belong to the same partition as $v$.
All edges that are added or deleted have one endpoint in $A$ and one in $B$,
and it follows that bipartiteness is preserved.\qed
\end{proof}

\begin{proposition}\label{prop:pivotbipartite}
Let $G$ be a bipartite graph, and let $\{u,v\} \in E$. Then $G^{(uv)}$ can be obtained
by ``toggling'' all edges between the sets $N_u \setminus \{v\}$ and 
$N_v \setminus \{u\}$, followed by a swapping of vertices $u$ and $v$.
\end{proposition}

\subsection{Coding Theory}\label{sec:coding}

A binary linear code, $\mathcal{C}$, is a linear subspace of $\GF(2)^n$ of dimension
$k$, where $0 \le k \le n$. $\mathcal{C}$~is called an $[n,k]$ code, and
the $2^k$ elements of $\mathcal{C}$ are called \emph{codewords}.
The \emph{Hamming weight} of $\boldsymbol{u} \in \GF(2)^n$, denoted $\wt(\boldsymbol{u})$,
is the number of nonzero components of $\boldsymbol{u}$.
The \emph{Hamming distance} between $\boldsymbol{u}, \boldsymbol{v} \in \GF(2)^n$
is $\wt(\boldsymbol{u} - \boldsymbol{v})$.
The \emph{minimum distance} of the code $\mathcal{C}$ is the minimal Hamming distance
between any two codewords of $\mathcal{C}$. Since $\mathcal{C}$ is a linear code,
the minimum distance is also given by the smallest weight of any codeword in $\mathcal{C}$.
A code with minimum distance~$d$ is called an $[n,k,d]$ code.
A code is \emph{decomposable} if it can be written as the \emph{direct sum} of two smaller codes.
For example, let $\mathcal{C}$ be an $[n,k,d]$ code and $\mathcal{C}'$ an $[n',k',d']$ code. The
direct sum, $\mathcal{C} \oplus \mathcal{C}' = \{u||v \mid u \in \mathcal{C}, v \in \mathcal{C}'\}$,
where $||$ means concatenation, is an $[n+n',k+k',\min\{d,d'\}]$ code.
Two codes, $\mathcal{C}$ and $\mathcal{C}'$, are considered to be \emph{equivalent} if
one can be obtained from the other by some permutation of the coordinates, or equivalently,
a permutation of the columns of a generator matrix.
We define the \emph{dual} of the code $\mathcal{C}$ with respect to the
standard inner product, $\mathcal{C}^\perp = 
\{ \boldsymbol{u} \in \GF(2)^n \mid \boldsymbol{u} \cdot \boldsymbol{c}=0 
\text{ for all } \boldsymbol{c} \in \mathcal{C} \}$.
$\mathcal{C}$~is called \emph{self-dual} if $\mathcal{C} = \mathcal{C}^\perp$,
and \emph{isodual} if $\mathcal{C}$ is equivalent to $\mathcal{C}^\perp$.
Self-dual and isodual codes must have even length $n$, and dimension $k=\frac{n}{2}$.
The code $\mathcal{C}$ can be defined by a $k \times n$ \emph{generator matrix}, $C$, whose 
rows span $\mathcal{C}$. 
A set of $k$ linearly independent columns of $C$ is called an \emph{information set} of $\mathcal{C}$. 
We can permute
the columns of $C$ such that an information set makes up the first $k$ columns. 
By elementary row operations, this matrix can then be transformed into a matrix of the form 
$C' = (I \mid P)$, where $I$ is a $k \times k$ identity matrix, 
and $P$ is some $k \times (n-k)$ matrix. The matrix $C'$, which is said to be of \emph{standard form},
generates a code $\mathcal{C}'$ which is equivalent to $\mathcal{C}$. Every code is equivalent 
to a code with a generator matrix of standard form. The matrix $H' = (P^\text{T} \mid I)$, 
where $I$ is an $(n-k) \times (n-k)$ identity matrix is called the \emph{parity check matrix} of 
$\mathcal{C}'$. Observe that $G'{H'}^\text{T} = \boldsymbol{0}$, where $\boldsymbol{0}$ is the 
all-zero vector. It follows that $H'$ must be the generator matrix of $\mathcal{C}'^\perp$.

\section{ELC and Code Equivalence}\label{sec:equiv}

As mentioned earlier, LC orbits of graphs correspond to equivalence classes of 
self-dual quantum codes. We have previously classified all such codes of 
length up to~12~\cite{selfdual}, by classifying LC orbits of simple undirected graphs.
In this paper, we show that ELC orbits of bipartite graphs correspond to
the equivalence classes of binary linear codes. First we explain how a
binary linear code can be represented by a graph.

\begin{definition}[\hspace{1pt}\hspace{-1pt}\cite{curtis,rijmen}]\label{def:code}
Let $\mathcal{C}$ be a binary linear $[n,k]$ code with generator matrix $C = (I \mid P)$.
Then the code $\mathcal{C}$ corresponds to the $(k,n-k)$-bipartite graph on $n$ vertices with
adjacency matrix
\[
\Gamma = \begin{pmatrix}\boldsymbol{0}_{k\times k} & P \\ P^\text{T} & \boldsymbol{0}_{(n-k)\times (n-k)}\end{pmatrix},
\]
where $\boldsymbol{0}$ denote all-zero matrices of the specified dimensions.
\end{definition}

\begin{theorem}\label{thm:pivotswap}
Let $G=(V,E)$ be the $(k,n-k)$-bipartite graph derived from a standard form generator matrix
$C=(I \mid P)$ of the $[n,k]$ code $\mathcal{C}$.
Let $G'$ be the graph obtained by performing ELC on the edge $\{u,v\} \in E$,
followed by a swapping of vertices $u$ and $v$.
Then the code $\mathcal{C}'$ generated by $C'=(I \mid P')$ corresponding to $G'$
is equivalent to $\mathcal{C}$,
and can be obtained by interchanging coordinates $u$ and $v$ of $\mathcal{C}$.
\end{theorem}
\begin{proof}
Assume, without loss of generality, that $u \le k$ and $v > k$.
$C'$ can be obtained from $C$ by adding row $u$ to all rows in $N_v \setminus \{u\}$
and then swapping columns $u$ and $v$, where $N_v$ denotes the neighbourhood of $v$ in $G$.
These operations preserve the equivalence of linear codes.
As described in Proposition~\ref{prop:pivotbipartite},
the bipartite graph $G$ is transformed into $G'$ by ``toggling'' all pairs of vertices 
$\{x,y\}$, where $x \in N_u \setminus \{v\}$ and $y \in N_v \setminus \{u\}$.
This action on the submatrix $P$ is implemented by the row additions on $C$ described above.
However, this also ``toggles'' the pairs $\{v,y\}$, where $y \in N_v \setminus \{u\}$,
transforming column $v$ of $C$ into a vector with 0 in all coordinates except $u$.
But column $u$ of $C$ now contains the original column $v$, and thus swapping
columns $u$ and $v$ restores the neighbourhood of $v$, giving the desired submatrix $P$.\qed
\end{proof}

\begin{corollary}
Applying any sequence of ELC operations to a graph $G$ corresponding to a code $\mathcal{C}$
will produce a graph corresponding to a code equivalent to $\mathcal{C}$.
\end{corollary}

Instead of mapping the generator matrix $C = (I\mid P)$ to the adjacency matrix of a 
bipartite graph in order to perform ELC on the edge $\{u,v\}$, we can work directly with 
the submatrix $P$.
Let the rows of $P$ be labeled $1,2,\ldots,k$ and the columns of $P$ be labeled $k+1,k+2,\ldots,n$.
Assume that $u$ indicates a row of $P$ and that $v$ indicates a column of $P$. 
The element $P_{ij}$ is then replaced by $1-P_{ij}$ if $i \ne u$, $j \ne v$, and $P_{uj} = P_{iv} = 1$.

\begin{figure}
 \centering
 {\hfill
 \subfloat[The Graph $G$]{\includegraphics[width=.30\linewidth]{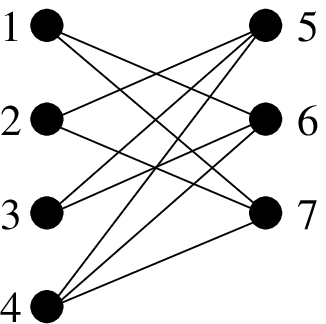}
   \label{fig:hamming1}}
 \hfill
 \subfloat[The Graph $G^{(27)}$]{\includegraphics[width=.30\linewidth]{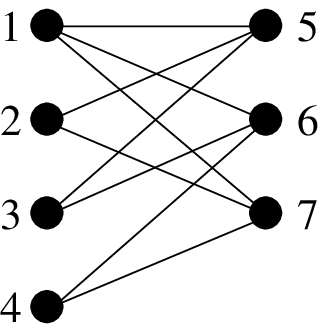}
   \label{fig:hamming2}}
 \hfill}
 \caption{Two Graph Representations of the $[7,4,3]$ Hamming Code}\label{fig:hammings}
\end{figure}

\begin{example}\label{ex:hamming}
The $[7,4,3]$ Hamming code has a generator matrix
\[
C = \left(
\begin{array}{cccc|ccc}
1&0&0&0& 0&1&1\\
0&1&0&0& 1&0&1\\
0&0&1&0& 1&1&0\\
0&0&0&1& 1&1&1
\end{array}
\right),
\]
which corresponds to the graph shown in Fig.~\ref{fig:hamming1}.
ELC on the edge $\{2,7\}$ produces the graph shown in Fig.~\ref{fig:hamming2}, which corresponds to the generator matrix
\[
C' = \left(
\begin{array}{cccc|ccc}
1&0&0&0& 1&1&1\\
0&1&0&0& 1&0&1\\
0&0&1&0& 1&1&0\\
0&0&0&1& 0&1&1
\end{array}
\right).
\]
The code generated by $C'$ is also obtained by swapping coordinates 2 and 7 of the code generated by $C$.
\end{example}

Consider a code $\mathcal{C}$.  As described in Section~\ref{sec:coding}, 
it is possible to go from a generator matrix of standard form, $C = (I \mid P)$, to
another generator matrix of standard form, $C'$, of a code equivalent to $\mathcal{C}$
by one of the $n!$ possible permutations of the columns of $C$, followed by
elementary row operations.
More precisely, we can get from $C$ to $C'$ via a combination of the following operations.
\begin{enumerate}
\item Permuting the columns of $P$.
\item Permuting the columns of $I$, followed by the same permutation on the rows of $C$, 
      to restore standard form.
\item Swapping columns from $I$ with columns from $P$, such that the first $k$ columns still
      is an information set, followed by some elementary row operations to restore standard form.
\end{enumerate}

\begin{theorem}\label{thm:alwayspivot}
Let $\mathcal{C}$ and $\mathcal{C}'$ be equivalent codes. Let $C$ and $C'$ be matrices
of standard form generating $\mathcal{C}$ and $\mathcal{C}'$.
Let $G$ and $G'$ be the bipartite graphs corresponding to $C$ and $C'$.
$G'$ is isomorphic to a graph obtained by performing some sequence of ELC operations on $G$.
\end{theorem}
\begin{proof}
$\mathcal{C}$ and $\mathcal{C}'$ must be related by a combination of the operations 1, 2, and 3 listed above.
It is easy to see that operations 1 and 2 applied to $G$ produce an isomorphic graph. It remains to 
prove that operation 3 always corresponds to some sequence of ELC operations.
We know from Theorem~\ref{thm:pivotswap} that swapping columns $u$ and $v$ of $C$, where
$u$ is part of $I$ and $v$ is part of $P$, corresponds to ELC on the edge $\{u,v\}$ of $G$,
followed by a swapping of the vertices $u$ and $v$.
When $\{u,v\}$ is not an edge of $G$, we can not swap columns $u$ and $v$ of $C$ via ELC.
In this case, coordinate $v$ of column $u$ is 0, and column $u$ has 1 in coordinate $u$ and 0 elsewhere.
Swapping these columns would result in a generator matrix where the first $k$ columns all have 0 
at coordinate $u$. These columns can not correspond to an information set.
It follows that if $\{u,v\}$ is not an edge of $G$, swapping columns $u$ and $v$ is not a valid
operation of type 3 in the above list. Thus ELC and graph isomorphism cover all possible operations
that map standard form generator matrices of equivalent codes to each other.\qed
\end{proof}

Let us for a moment consider ELC orbits of \emph{labeled} graphs, i.e., 
where we do not take isomorphism into consideration. 
Let $G=(V,E)$ be the connected bipartite graph representing the indecomposable code $\mathcal{C}$,
and $G^{(uv)}$ be the graph obtained by ELC on the edge $\{u,v\} \in E$.
Since we perform ELC on $\{u,v\}$ without swapping $u$ and $v$ afterwards,
the adjacency matrix of $G^{(uv)}$ will not be of the type we saw in 
Definition~\ref{def:code}. Assuming that vertices $\{1,2,\ldots,k\}$ make up one of the partitions
of the bipartite graph $G$, we can think of $G$ as a graph corresponding to the information set
$\{1,2,\ldots,k\}$ of $\mathcal{C}$. Assume that $u \le k$ and $v > k$.
$G^{(uv)}$ will then represent another information set of $\mathcal{C}$, namely 
$\{1,2,\ldots,k\} \setminus \{u\} \cup \{v\}$.

\begin{theorem}\label{thm:infosets}
Let $G$ be a connected bipartite graph representing the indecomposable code $\mathcal{C}$.
Each labeled graph in the ELC orbit of $G$ corresponds to an information set of $\mathcal{C}$.
If $\mathcal{C}$ is a self-dual code, each graph corresponds to two information sets, one for
each partition.
Moreover, the number of information sets of $\mathcal{C}$ equals the number of
labeled graphs in the ELC orbit of $G$, or twice the number of graphs if $\mathcal{C}$ is a self-dual code.
\end{theorem}
\begin{proof}
Performing ELC without swapping vertices afterwards corresponds to elementary row operations on the associated
generator matrix, and will thus leave the code invariant. The only thing we change with ELC is the information
set of the code, as indicated by the bipartition of the graph.
We know from Theorem~\ref{thm:alwayspivot} that if two generator matrices of standard form generate equivalent
codes, we can always get from one to the other via ELC operations on the associated graph.
It follows from this that when we consider labeled graphs, and do not swap vertices to obtain a code
of standard form, we find all information sets in the ELC orbit.
If and only if a code is self-dual, $(I\mid P)$ will generate the same code as $(P^\text{T} \mid I)$.
Since the matrices $(I\mid P)$ and $(P^T \mid I)$ correspond to exactly the same graph, but two different information
sets, we must multiply the ELC orbit size with two to get the number of information sets of a self-dual code.\qed
\end{proof}

Note that the distinction between ELC with or without a final swapping of vertices is only
significant when we want to find information sets. For other applications, where we consider
graphs up to isomorphism, this distinction is not of importance.

\begin{theorem}
The minimum distance, $d$, of a binary linear $[n,k,d]$ code $\mathcal{C}$, is equal to $\delta + 1$, where
$\delta$ is the smallest vertex degree of any vertex in the partition of size $k$
over all graphs in the associated ELC orbit.
\end{theorem}
\begin{proof}
If there is a vertex with degree~$d-1$, belonging to the partition of size~$k$, in the ELC orbit, there
is a row of weight~$d$ in a generator matrix that generates a code equivalent to $\mathcal{C}$.
Hence there must also be a codeword of weight $d$ in $\mathcal{C}$.
We need to show that when $d$ is the minimum distance of $\mathcal{C}$,
such a vertex always exists. 
Let $C$ be the standard form generator matrix of $\mathcal{C}$. If $C$ contains a row
of weight $d$, we are done. Otherwise, select a codeword $\boldsymbol{c}$ of weight $d$, generated by $C$, 
and let the $i$-th row of $C$ be one of the rows that $\boldsymbol{c}$ is linearly dependent on.
Replace the $i$-th row of $C$ by $\boldsymbol{c}$ to get $C'$. 
Permute the columns of $C'$ to obtain $C''$ where the first $k$ columns is still an information set,
and where $\boldsymbol{c}$ is mapped to $\boldsymbol{c}'$ with 1 in coordinate $i$, 
with the rest of the $k$ first coordinates being 0.
That such a permutation will always exist follows from the fact that
$\boldsymbol{c}$ has weight~$d$ while all other rows of $C'$ have weight greater than $d-1$ in the 
last $n-k$ coordinates.
Thus, for each coordinate $j \le k$, $j \ne i$, where $\boldsymbol{c}$ 
is~1, there must exist a distinct coordinate $l > k$ where $\boldsymbol{c}$ is~0 and 
the $j$-th row of $C'$ is~1.
We can transform $C''$ into a matrix of the form $(I \mid P)$ by elementary row operations.
Row $i$ of this final matrix has weight~$d$, and thus the corresponding bipartite graph
has a vertex with degree $d-1$.\qed
\end{proof}

\section{Classification of ELC Orbits}\label{sec:class}

We have previously classified all self-dual additive codes over $\GF(4)$ of length up 
to~12~\cite{selfdual,mscthesis}, by classifying orbits of simple undirected graphs
with respect to local complementation and graph isomorphism.
In Table~\ref{tab:lcorbits}, the sequence $(i_n^{LC})$ gives the number of LC orbits 
of connected graphs on $n$ vertices, while $(t_n^{LC})$ gives the total number of 
LC orbits of graphs on $n$ vertices.
A database containing one representative from each LC orbit is available 
at \url{http://www.ii.uib.no/~larsed/vncorbits/}.

\begin{table}
\centering
\small{
\centering
\caption{Numbers of LC Orbits}
\label{tab:lcorbits}
\begin{tabular}{ccccccccccccc}
\toprule
$n$ & 1 & 2 & 3 & 4 & 5 & 6 & 7 & 8 & 9 & 10 & 11 & 12 \\
\midrule
$i_n^{LC}$ & 1 & 1 & 1 & 2 &  4 & 11 & 26 & 101 & 440 & 3,132 & 40,457 & 1,274,068 \\
$t_n^{LC}$ & 1 & 2 & 3 & 6 & 11 & 26 & 59 & 182 & 675 & 3,990 & 45,144 & 1,323,363 \\
\bottomrule
\end{tabular}
}
\end{table}

By recursively applying ELC operations 
to all edges of a graph, whilst checking for graph isomorphism using the program 
\emph{nauty}~\cite{nauty}, we can find all members of the ELC orbit. 
Let $\boldsymbol{G}_n$ be the set of all unlabeled simple undirected connected 
graphs on $n$ vertices. Let the set of all distinct ELC orbits of connected graphs 
on $n$ vertices be a partitioning of $\boldsymbol{G}_n$ into $i_n^{ELC}$ disjoint sets.
Our previous classification of the LC orbits of all graphs of up to 12 vertices helps 
us to classify ELC orbits, since it follows from Definition~\ref{prop:triplelc}
that each LC orbit can be partitioned into a set of disjoint ELC orbits. 
We have used this fact to classify all ELC orbits of graphs on up to 12 vertices,
a computation that required approximately one month of running time on a parallel
cluster computer.
In Table~\ref{tab:pivorbits}, the sequence $(i_n^{ELC})$ gives the number of ELC orbits 
of connected graphs on $n$ vertices, while $(t_n^{ELC})$ gives the total number of 
ELC orbits of graphs on $n$ vertices.
Note that the value of $t_n$ can be derived easily once the sequence $(i_m)$ is known 
for $1 \le m \le n$, using the \emph{Euler transform}~\cite{sloane2},
\begin{eqnarray*}
c_n &=& \sum_{d|n} d i_d,\\
t_1 &=& c_1,\\
t_n &=& \frac{1}{n}\left( c_n + \sum_{k=1}^{n-1} c_k t_{n-k} \right).
\end{eqnarray*}
A database containing one representative from each ELC orbit
can be found at \url{http://www.ii.uib.no/~larsed/pivot/}.

\begin{table}
\centering
{\small
\centering
\caption{Numbers of ELC Orbits and Binary Linear Codes}
\label{tab:pivorbits}
\begin{tabular}{rrrrrrr}
\toprule
$n$ & $i_n^{ELC}$    & $t_n^{ELC}$      & $i_n^{ELC,B}$& $t_n^{ELC,B}$ & $i_n^{C}$ & $i_n^{C_{iso}}$\\
\midrule
 1 &           1 &           1 &          1 &         1 &          1 &      - \\
 2 &           1 &           2 &          1 &         2 &          1 &      1 \\
 3 &           2 &           4 &          1 &         3 &          2 &      - \\
 4 &           4 &           9 &          2 &         6 &          3 &      1 \\
 5 &          10 &          21 &          3 &        10 &          6 &      - \\
 6 &          35 &          64 &          8 &        22 &         13 &      3 \\
 7 &         134 &         218 &         15 &        43 &         30 &      - \\
 8 &         777 &       1,068 &         43 &       104 &         76 &     10 \\
 9 &       6,702 &       8,038 &        110 &       250 &        220 &      - \\
10 &     104,825 &     114,188 &        370 &       720 &        700 &     40 \\
11 &   3,370,317 &   3,493,965 &      1,260 &     2,229 &      2,520 &      - \\
12 & 231,557,290 & 235,176,097 &      5,366 &     8,361 &     10,503 &    229 \\
13 &   \emph{?}  &    \emph{?} &     25,684 &    36,441 &     51,368 &      - \\
14 &             &             &    154,104 &   199,610 &    306,328 &  1,880 \\
15 &             &             &  1,156,716 & 1,395,326 &  2,313,432 &      - \\
16 &             &             &   \emph{?} &  \emph{?} & \emph{23,069,977} & \emph{?} \\
17 &             &             &   \emph{157,302,628}  &  \emph{?}  & \emph{314,605,256}    &      - \\
\bottomrule
\end{tabular}
}
\end{table}

We are particularly interested in bipartite graphs, because of their connection to binary linear codes.
For the classification of the orbits of bipartite graphs with respect to ELC and graph
isomorphism, the following technique is helpful.
If $G$ is an $(a,b)$-bipartite graph, it has $2^a + 2^b -2$ possible \emph{extensions}.
Each extension is formed by adding a new vertex and joining it to all possible combinations 
of at least one of the old vertices.
Let $\boldsymbol{P}_n$ be a set containing one representative from each
ELC orbit of all connected bipartite graphs on $n$ vertices. The set $\boldsymbol{E}_{n}$ 
is formed by making all possible extensions of all graphs in $\boldsymbol{P}_{n-1}$.
It can then be shown that $\boldsymbol{P}_n \subset \boldsymbol{E}_{n}$, i.e., that
the set $\boldsymbol{E}_{n}$ will contain at least one representative from each
ELC orbit of connected bipartite graphs on $n$ vertices.
The set $\boldsymbol{E}_{n}$ will be much smaller than $\boldsymbol{G}_n$, so it will be more
efficient to search for a set of ELC orbit representatives within $\boldsymbol{E}_{n}$.
A similar technique was used by Glynn, et~al.~\cite{glynnbook} to classify LC orbits.

In Table~\ref{tab:pivorbits}, the sequence $(i_n^{ELC,B})$ gives the number of 
ELC orbits of connected bipartite graphs on $n$ vertices, and
$(t_n^{ELC,B})$ gives the total number of ELC orbits of bipartite graphs on $n$ vertices.
A database containing one representative from each of these orbits
can be found at \url{http://www.ii.uib.no/~larsed/pivot/}.

\begin{theorem}\label{prop:codenumbers}
Let $k \ne \frac{n}{2}$. Then the number of inequivalent binary linear $[n,k]$ codes, 
which is also the number of inequivalent $[n,n-k]$ codes, is equal to the number of
ELC orbits of $(n-k,k)$-bipartite graphs.

When $n$ is even and $k = \frac{n}{2}$, the number of inequivalent binary linear $[n,k]$ codes
is equal to twice the number of ELC orbits of $(k,k)$-bipartite graphs minus
the number of isodual codes of length $n$.
\end{theorem}
\begin{proof}
We recall that if a code $\mathcal{C}$ is generated by $(I\mid P)$, then its dual,
$\mathcal{C}^\perp$, is generated by $(P^\text{T} \mid I)$. Also note that $\mathcal{C}^\perp$ is
equivalent to the code generated by $(I \mid P^\text{T})$.
The bipartite graphs corresponding to the codes generated by $(I\mid P)$ and $(I\mid P^\text{T})$ 
are isomorphic. It follows that the ELC orbit associated with an $[n,k]$ code $\mathcal{C}$ is 
simultaneously the orbit associated with the dual $[n,n-k]$ code $\mathcal{C}^\perp$. 
In the case where $k = \frac{n}{2}$, each ELC orbit corresponds to two non-equivalent $[n,k]$ codes,
except in the case where $\mathcal{C}$ is isodual.\qed
\end{proof}

\begin{corollary}
The total number of binary linear codes of length $n$ is equal to twice the number of
ELC orbits of bipartite graphs on $n$ vertices, minus the number of isodual codes of length $n$.
\end{corollary}

Note that if we only consider connected graphs on $n$ vertices,
we get the number of indecomposable codes of length $n$, $i_n^{C}$, i.e.,
the codes that can not be written as the direct sum of two smaller codes.
The total number of codes can easily be 
derived from the values of $(i_n^{C})$. Table~\ref{tab:pivorbits}
gives the number of ELC orbits of connected bipartite graphs on $n$ vertices, $i_n^{ELC,B}$, the
number of indecomposable binary linear codes of length $n$, $i_n^{C}$, and
the number of indecomposable isodual codes of length $n$, $i_n^{C_{iso}}$.
A method for counting the number of binary linear codes by using computer algebra tools
was devised by  Fripertinger and Kerber~\cite{isometry}. A table enumerating binary linear codes
of length up to 25 is available online at \url{http://www.mathe2.uni-bayreuth.de/frib/codes/tables_2.html}.
The numbers in italics in Table~\ref{tab:pivorbits} are taken from this webpage.
Note however that this approach only gives the number of inequivalent codes, and does
not produce the codes themselves. 
Classification of all binary linear codes of length up to 14 and with distance at least 3
was carried out by Östergård~\cite{ostergard}. He also used a graph-based algorithm, but
one quite different from the method described in this paper.
In a recent book by Kaski and Östergård~\cite{ostergard2}, it is proposed as a research problem to
extend this classification to lengths higher than 14.
Sang-il Oum [personal communication] demonstrated that the 1,395,326 ELC orbits of 
bipartite graphs on 15 vertices can be generated in about 58 hours.
This indicates that classification of codes by ELC orbits is comparable to the currently
best known algorithm. It may also be possible that our method will be more efficient than 
existing algorithms for classifying special types of codes. For instance, matrices
of the form $(I \mid P)$, where $P$ is symmetric, generate a subset of the isodual codes.
The bipartite graphs corresponding to these codes, which  were also studied by 
Curtis~\cite{curtis}, should be well suited to our method, since any graph of this type must arise
as an extension of a graph of the same type.

\paragraph*{Acknowledgements}
This research was supported by the Research Council of Norway.
We would like to thank the Bergen Center for Computational Science, whose
cluster computer made the results in this paper possible.
Thanks to Joakim G. Knudsen for help with improving Theorem~\ref{thm:infosets}.

\end{document}